\documentclass[a4paper,10pt]{amsart}
\usepackage{epsf}  
\usepackage{amsfonts, amssymb, amsthm, amsmath}  
\usepackage{hyperref}
\usepackage{tikz}
\usetikzlibrary{matrix,arrows,decorations.pathmorphing}
\usepackage{tikz-cd}
\usepackage{graphicx}
\usepackage{psfrag}
\newtheorem{thm}{Theorem}

\newtheorem{remark}[thm]{Remark}  
\newtheorem{defn}[thm]{Definition}

\newtheorem{example}[thm]{Example}  
\numberwithin{thm}{section}

\providecommand{\Y}{\mathcal Y}

\providecommand{\totl}[1]{\ensuremath{\lceil #1\rceil }}
\providecommand{\totb}[1]{\ensuremath{\underline{ #1}}}

\DeclareMathOperator{\rend}{End}
\DeclareMathOperator{\Aut}{Aut}

\newcommand{\rhf} {{}^{\phantom{f}r}_{fg}H}
\newcommand{\rh}{{}^{r}H}
\newcommand{\ie}{\text{ied}}
\newcommand{\ed}{\text{ed}}
\newcommand{\ex}{\bold}
\providecommand {\e}[1]{\mathfrak t^{#1}}
\providecommand{\C}[2]{\ensuremath {C^{#1,\underline{#2}}}}

\newcommand{\tc}[1]{\check\rvert_{#1}}

\newcommand{\Mod}{\mathcal M}

\newcommand{\Msw}{\mathcal M^{st}}

\DeclareMathOperator{\id}{id}
\DeclareMathOperator{\expl}{Expl}

\newcommand{\dbar}{\bar{\partial}}

\providecommand{\et}[2]{\ensuremath{\bold T^{#1}_{#2}}}
\providecommand{\lrb}[1]{\ensuremath{\left(#1\right)}}
\providecommand{\abs}[1]{\left\lvert #1\right\rvert}

\author{Brett Parker   }
\email{brettdparker@gmail.com}  
\thanks{ Partially funded by the Simons center, and partially funded by ARC grant  DP140100296.}
  
\title[Exploded manifolds and tropical gluing formula]{Notes  on Exploded manifolds and a tropical gluing formula for Gromov--Witten invariants}

\begin{document}
\maketitle

\begin{abstract}Notes for a short lecture series, covering exploded manifolds, the moduli stack of curves in exploded manifolds, and a tropical gluing formula for Gromov--Witten invariants --- a gluing formula providing a degeneration formula for Gromov--Witten invariants in normal-crossing degenerations. I gave the original lecture series in April 2016 at the Simons Center for Geometry and Physics at Stonybrook. Video of the lectures is available on the  \href{http://scgp.stonybrook.edu/video_portal/video.php?id=2595}{SCGP website}. 

\end{abstract}

\tableofcontents

\newpage

\section{Introduction}

If we want a gluing or degeneration formula for Gromov--Witten invariants in a normal-crossing degeneration, we must keep track of information on two different scales: a small scale involving differential geometry of manifolds with products of cylindrical ends, and a large scale involving piecewise integral-affine geometry. On this large scale, holomorphic curves look like tropical curves, and consequently the gluing formula for Gromov--Witten invariants is a sum over tropical curves. 
This large and small scale information appears simultaneously  in the degenerating family of targets, the domain of holomorphic curves, and the moduli stack of holomorphic curves, so it is natural and convenient to work in a category  that keeps track of this large and small scale information systematically: the category of exploded manifolds.

 Without using an appropriate category, such as exploded manifolds or log schemes, an explanation of Gromov--Witten invariants relative normal-crossing divisors or the associated gluing formula will seen surprising, unnatural, and unmotivated.\footnote{Some motivation is given in the  introduction to \cite{iec} or slides from talks  linked to from my website,  \href{http://maths-people.anu.edu.au/~parkerb}{http://maths-people.anu.edu.au/$\sim$parkerb} --- the actual lecture series started with an unrecorded talk similar to the one at  \href{https://prezi.com/yj74jcoy44vi/}{this link}. The gluing formula is also described in an easy case without special language in \cite{tpgf}. }    By contrast, once  exploded manifolds are understood, Gromov--Witten invariants relative normal-crossing divisors define themselves naturally with no creative leaps required\footnote{This is not to say that defining these Gromov Witten invariants is easy. In both the log and exploded settings, it requires hard, technical work. For log Gromov--Witten invariants, see \cite{GSlogGW, Chen, acgw}. A separate, but related approach to defining GW invariants relative a version of normal-crossing divisors in the symplectic setting is found in \cite{IonelGW}. For the relationship between these approaches and the approach using exploded manifolds, see \cite{elc,tropicalIonel}.}, and the associated tropical gluing formula is inevitable. We replace the original degenerating family of manifolds with a smooth family of exploded manifolds.   In this exploded family, holomorphic curves are holomorphic maps of curves within the category of exploded manifolds. The moduli stack of such holomorphic curves is a nice stack over the category of exploded manifolds, and is naturally smooth (modulo the usual transversality considerations). This moduli stack carries a natural virtual fundamental class compatible with pullbacks and therefore defining Gromov--Witten invariants that do not change in families. The Gromov--Witten invariants of any  fiber are the Gromov--Witten invariants of the original smooth manifold, and when we calculate these invariants for a `singular'  fiber, we find our tropical gluing formula.

In order to explain relative Gromov--Witten invariants and the tropical gluing formula, we first study exploded manifolds.
 Exploded manifolds may seem strange, but  their strangest features are forced by the geometry of holomorphic curves in normal-crossing degenerations. These geometric constraints created an unintentional link between exploded manifolds and log schemes. 
  
  These notes accompany my  series of 3 lectures recorded in April 2016 at the Simons Center for Geometry and Physics at Stonybrook. The videos are available on the \href{http://scgp.stonybrook.edu/video_portal/video.php?id=2595}{SCGP website}. Video of an introductory talk  I refer to in these recordings is not available, but slides from analogous talks are available at \href{https://prezi.com/yj74jcoy44vi/}{this link} or \href{https://prezi.com/fweevmpxkthu/ }{this link}. Please do not hesitate to contact me if you have any questions or comments.

\section{Exploded manifolds}

We begin with the definition of an exploded manifold, however we can extract little from it without further definitions and  examples. 
\begin{defn}An exploded manifold is an abstract exploded space locally isomorphic to $\mathbb R^{n}\times \et mP$.
\end{defn}
Obviously, we need to understand what an abstract exploded space is, and what  these $\mathbb R^{n}\times \et mP$ are. All that can be extracted from the above definition is that exploded manifolds  have coordinate charts isomorphic to $\mathbb R^{n}\times \et mP$. The type of these coordinate charts depends on nonnegative integers $n$ and $m$ and a $m$--dimensional integral-affine polytope $P$ (defined later). We shall see that the (real) dimension of such a coordinate chart is $n+2m$.

\begin{defn}An abstract exploded space  is \begin{itemize}\item a set $\ex B$ with topology induced from a surjective map to a Hausdorff topological space $\ex B\longrightarrow \totl{\ex B}$, \item and a sheaf, $\mathcal E^{\times}$, of $\mathbb C^{*}\e{\mathbb R}$--valued functions, containing constants, and closed under multiplication and taking inverses, where $\mathbb C^{*}\e{\mathbb R}$ means the group $(\mathbb C^{*},\times)\times(\mathbb R,+)$.\end{itemize} \end{defn}

A morphism  of abstract exploded spaces is a continuous map $\phi:\ex A\longrightarrow \ex B$ so that $\phi^{*}(\mathcal E^{\times}(\ex B))\subset \mathcal E^{\times}(\ex A)$. We call functions in $\mathcal E^{\times}$ exploded functions.

We use $\mathbb C^{*}\e{\mathbb R}$ for $(\mathbb C^{*},\times)\times(\mathbb R,+)$ because it is the group of units for a semiring $\mathbb C\e{\mathbb R}$ --- an algebraic  exploded manifold is like a nice algebraic space over this exploded semiring $\mathbb C\e{\mathbb R}$. As a set, $\mathbb C\e{\mathbb R}$ is  $\mathbb C\times \mathbb R$, with   $(c,a)$ written $c\e{a}$. The laws for multiplication and addition are then
\[c_{1}\e{a_{1}}\times c_{2}\e{a_{2}}:=(c_{1}c_{2})\e{a_{1}+a_{2}}\]
and
\[c_{1}\e{a_{1}}+c_{2}\e{a_{2}}:=\begin{cases}c_{1}\e{a_{1}}\ \ \ \ \ \ \ \ \ \text{ if }a_{1}<a_{2}
\\ (c_{1}+c_{2})\e{a_{1}}\text{ if }a_{1}=a_{2}
\\ c_{2}\e{a_{2}}\ \ \ \ \ \ \ \ \ \text{ if }a_{2}<a_{1}\end{cases}\]
The same construction works with any ring replacing $\mathbb C$ --- for example $0\e{\mathbb R}$ is the tropical semiring. There is a tropical-part homomorphism from $\mathbb C^{*}\e{\mathbb R}$ to the tropical numbers.
\[\totb{c\e a}:=a\]
Equally important, is the following smooth-part homomorphism only defined on $\mathbb C\e{[0,\infty)}$, the set of elements $c\e a$ where $a\geq 0$.
\[\totl{c\e a}:=\begin{cases}0\text{ if }a>0
\\ c\text{ if }a=0\end{cases}\]
The addition law on $\mathbb C^{*}\e {\mathbb R}$ and the smooth-part homomorphism make sense if we think $\e{}$ is infinitesimal.

There are analogous tropical and smooth-part maps for exploded manifolds. The surjective map to a Hausdorff topological space, $\ex B\longrightarrow \totl{\ex B}$, is related to the above smooth-part homomorphism, and $\totl{\ex B}$ is called the smooth part of $\ex B$. Corresponding to the tropical-part homomorphism, there is a surjective tropical-part map $\ex B\longrightarrow \totb{\ex B}$, where $\totb{\ex B}$ is some complex of integral-affine polytopes called the tropical part of $\ex B$. The tropical part of  a coordinate chart $\mathbb R^{n}\times \et mP$ is $P$, so coordinate charts on exploded manifolds are classified by dimension and their tropical part.

We now examine  examples of exploded manifolds, building towards understanding a general coordinate chart $\mathbb R^{n}\times \et mP$.

\begin{example}$\ex T:=\et 1{\mathbb R}$
\end{example}
As a topological space, $\ex T=\mathbb C^{*}\e{\mathbb R}$ with the trivial indiscrete topology. Let $z:\ex T\longrightarrow \mathbb C^{*}\e{\mathbb R}$ be the corresponding coordinate. The exploded functions on $\ex T$ are the monomials in $z$.
\[\mathcal E^{\times}(\ex T):=\{cz^{n}\text{ so that }n\in \mathbb Z \text{ and }c\in \mathbb C^{*}\e{\mathbb R}\}\]
This example is both trivial and important: on any abstract exploded space, the sheaf of exploded functions is the sheaf of morphisms to  $ \ex T$.

\begin{example}$\et 1{[0,\infty)}:=\expl (\mathbb C, 0)$
\end{example}
As a set,  $\et 1{[0,\infty)}:=\mathbb C^{*}\e{[0,\infty)}$, with topology induced from the smooth-part homomorphism $\mathbb C^{*}\e{[0,\infty)}\longrightarrow\mathbb C:=\totl{\et 1{[0,\infty)}}$. Let $z:\et 1{[0,\infty)}\longrightarrow \mathbb C^{*}\e{[0,\infty)}\subset \mathbb C^{*}\e{\mathbb R}$ be the corresponding coordinate on $\et 1{[0,\infty)}$, and let $\totl z:\et 1{[0,\infty)}\longrightarrow \mathbb C$ be the composition of $z$ with the smooth-part homomorphism. The exploded functions on $\et 1{[0,\infty)}$ are all functions in the form
\[h(\totl{z})\e az^{n}\]
where $h$ is a smooth\footnote{We can put a different regularity on $\et 1{[0,\infty)}$, such as holomorphic, or continuous, by changing the regularity required of $h$. Holomorphic exploded manifolds are abstract exploded spaces locally isomorphic to \emph{an open subset of} $\mathbb C^{n}\times \et {}{P}$ with its sheaf of holomorphic exploded functions. Later on, we shall need to use a regularity --- weaker, but pretty much as good as smooth --- called $\C\infty1$. See Remark \ref{C}.}, $\mathbb C^{*}$--valued function, $a\in\mathbb R$ is locally constant, and $n\in\mathbb Z$.   

Later on, we discuss the explosion functor --- a functor producing a holomorphic  exploded manifold from a complex manifold with normal-crossing divisors ( or more generally, a complex,  log smooth, log scheme). The explosion of $\mathbb C$ with the divisor $0$ is $\et 1{[0,\infty)}$. Think of this explosion as `exploding' the divisor $0$ and replacing it with $\mathbb C^{*}\e{(0,\infty)}$,  replacing the divisor with  a collection of cylinders, $\mathbb C^{*}$, indexed by $(0,\infty)$.

\begin{remark}\label{log remark}[Relationship of exploded manifolds to log schemes]
For exploded manifolds $\ex B$,  the sheaf of morphisms to $\et1{[0,\infty)}$ is the sheaf of exploded functions with tropical part in $[0,\infty)\subset\mathbb R$. We can compose such exploded functions with the smooth-part homomorphism to obtain $\mathbb C$--valued functions. Such $\mathbb C$--valued functions come from functions on $\totl{\ex B}$. So, maps from $\ex B$ to $\et1{[0,\infty)}$ form a sheaf of monoids on $\totl{\ex B}$ with a homomorphism to the sheaf of $\mathbb C$--valued functions on $\totl{\ex B}$. Such a structure is called a log structure. Exploded manifolds (with tropical part consisting of polytopes not containing any entire lines) may be defined using this log structure together with the natural map to a point $p^{\dag}$ with log structure  the smooth-part homomorphism $\mathbb C^{*}\e{[0,\infty)}\longrightarrow\mathbb C$. From this perspective, the explosion of a complex  log scheme $M$ is $M\times p^{\dag}$ with the natural projection to $p^{\dag}$. For further details on the relationship between exploded manifolds and log schemes, see \cite{elc}.
\end{remark}

After this warmup, we are ready to define $\et mP$. Let $P\subset \mathbb R^{m}$ be an integral-affine polytope --- a subset of $\mathbb R^{m}$ with nonempty interior, cut out by finitely many integral-affine inequalities in the form
\[x\cdot\alpha+a\geq 0 \text{ or } x\cdot\alpha+a>0\]
where  $\alpha\in\mathbb Z^{m}$ and $a\in\mathbb R$.

\begin{example}$\et mP:=\et{}P$\end{example}
As a set, $\et mP$ is the subset of $\lrb{\mathbb C^{*}\e{\mathbb R}}^{m}$ with tropical part in $P$. 
\[\et mP:=\left\{(c_{1}\e {a_{1}},\dotsc, c_{m}\e{a_{n}})\in \lrb{\mathbb C^{*}\e{\mathbb R}}^{m} \text{ so that }(a_{1},\dotsc,a_{m})\in P\right\}\]
Accordingly, we get coordinate functions $z_{i}:\et mP\longrightarrow \ex T$ so that $(\totb{z_{1}},\dotsc,\totb{z_{m}})$ defines our tropical-part map, $\et {}P\longrightarrow \totb{\et {}P}:=P$. 
Given an integral-affine map $P\longrightarrow[0,\infty)$ in the form $x\mapsto a+x\cdot \alpha$, there is a corresponding monomial on $\et mP$
\[\e az^{\alpha}:\et mP\longrightarrow \mathbb C^{*}\e{[0,\infty)}\ \ \ \ \ \ \ \ \text{where }z^{\alpha}:=\prod_{i=1}^{m}z_{i}^{\alpha_{i}}\ .\]
We can compose such a monomial with the smooth-part homomorphism to define a $\mathbb C$--valued function  called a smooth monomial.
\[\zeta:=\totl{\e az^{\alpha}}:\et {}P\longrightarrow \mathbb C\]
The definition of the smooth-part homomorphism implies that such a monomial is nonzero exactly where $\totb{\e az^{\alpha}}=0$. Such smooth monomials form a finitely generated monoid; choose some basis $\zeta_{1},\dotsc,\zeta_{n}$,   then, the smooth part of $\et {}P$ is  the image of $(\zeta_{1},\dotsc,\zeta_{n})$ in $\mathbb C^{n}$, and the smooth-part map $\et {}P\longrightarrow \totl{\et {}P}$ is $(\zeta_{1},\dotsc,\zeta_{n})$. The induced topology on $\et {}P$ is the coarsest topology with all smooth monomials continuous.

We now define the sheaf of exploded functions  $\mathcal E^{\times}(\et {}P)$. The exploded functions are all functions in the form 
\[h(\zeta_{1},\dotsc,\zeta_{n})\e a z^{\alpha}\]
where $h$ is smooth and $\mathbb C^{*}$--valued, and $a\in\mathbb R$ and $\alpha\in\mathbb Z^{m}$ are locally constant. Again, we can change the regularity of our exploded manifold by changing the regularity required of $h$. 

\

Once $\et {}P$ is defined, there are no surprises in defining a general coordinate chart $\mathbb R^{n}\times \et {}P$. As a topological space, this is the product of $\et {}P$ with $\mathbb R^{n}$, its smooth part is the product of $\totl{\et {}P}$ with $\mathbb R^{n}$, and its tropical part is $P$. The exploded functions on $\mathbb R^{n}\times \et {}P$ are as above, but now $h$ may also depend on $\mathbb R^{n}$.

\begin{remark}The reader should now be able to verify the following facts.

\begin{enumerate}
\item A map from any open subset $U$ of a coordinate chart to $\et mP$ is equivalent to $m$ exploded functions $f_{i}\in \mathcal E^{\times}(U)$ so that $(\totb {f_{1}},\dotsc,\totb{f_{m}})\in P$.
\item A map  $\et {}P\longrightarrow \et {}Q$ induces an integral-affine map $P\longrightarrow Q$. 
\item Given any integral-affine map $P\longrightarrow Q$ there exists a map $\et {}P\longrightarrow \et{}Q$ inducing it.
\item $\et{}P$ is isomorphic to $\et {}Q$ if and only if $P$ is isomorphic to $Q$ as an integral-affine polytope.
\item Any map $\et {}P\longrightarrow \et{}Q$ induces a map $\totl{\et {}P}\longrightarrow \totl{\et{}Q}$ compatible with the smooth structures on $\totl{\et {}P}$ and $\totl{\et {}Q}$.
\item For any exploded manifold $\ex B$, there exist functorial smooth-part and tropical-part maps, $\ex B\longrightarrow \totl{\ex B}$ and $\ex B\longrightarrow \totb{\ex B}$, that are as  described on coordinate charts. Any morphism  $\phi:\ex A\longrightarrow \ex B$ has a smooth part $\totl\phi$ and tropical part $\totb\phi$ so that the following diagram commutes.
\[\begin{tikzcd}\totl{\ex A}\rar{\totl{\phi}}&\totl{\ex B}
\\ \ex A\uar\dar\rar{\phi}&\ex B\uar\dar
\\ \totb{\ex A}\rar{\totb\phi}&\totb{\ex B}\end{tikzcd}\]
The tropical part can be complicated by how the tropical parts of different coordinate charts glue together. See section 4 of \cite{iec} for more details.
\end{enumerate}
\end{remark}

\begin{example}\label{et[0,l]}$\et{}{[0,l]}$\end{example}
As a set $\et{}{[0,l]}=\mathbb C^{*}\e{[0,l]}$. A basis for the smooth monomials on $\et{}{[0,l]}$ is $\zeta_{1}:=\totl{z}$ and $\zeta_{2}=\totl{\e lz^{-1}}$. The smooth part of $\et{}{[0,l]}$ is 
\[\totl{\et{}{[0,l]}}=\{\zeta_{1}\zeta_{2}=0\}\subset\mathbb C^{2}\ ,\]
two complex planes joined at the origin --- the local model around a node of a holomorphic curve (in the smooth category). Any holomorphic curve $\ex C$ in the category of exploded manifolds  has  smooth part  a nodal curve $\totl{\ex C}$. The nodes of $\totl{\ex C}$ correspond to strata of $\ex C$ we shall call internal edges. And,  around an internal edge, $\ex C$  is locally modeled on $\et {}{[0,l]}$. So,  over the node in $\totl{\ex C}$ is an entire stratum of $\ex C$ isomorphic to $\et{}{(0,l)}$. Because $\et{}{(0,l)}$ has no nonconstant smooth monomials,  the restriction of exploded functions to an edge are just the monomials $cz^{n}$, so maps of $\et{}{[0,1]}$ are rigid restricted to $\et{}{(0,l)}$. A key observation (which shall lead to our  gluing formula) is that choosing a map from $\et{}{[0,l]}$ is equivalent to choosing a maps from $\et{}{[0,l)}$ and $\et{}{(0,l]}$ agreeing  on $\et{}{(0,l)}$.

An instructive way to construct $\et{}{[0,l]}$ is as the subset of $\et{2}{[0,\infty)^{2}}$ where $z_{1}z_{2}=1\e l$, so $\zeta_{i}$ is the restriction of $\totl{z_{i}}$. This subset is a fiber of the map
\[z_{1}z_{2}:\et {}{[0,\infty)^{2}}\longrightarrow \et{}{[0,\infty)}\ .\]
The smooth part of this map is the local model for node formation.
\[\zeta_{1}\zeta_{2}:\mathbb C^{2}\longrightarrow \mathbb C\]
Over points $\epsilon\e{0}\in \et {}{[0,\infty)}$ we get smooth manifolds $\zeta_{1}\zeta_{2}=\epsilon$, but over points $c\e l$ where $l>0$, we get exploded manifolds isomorphic to $\et{}{[0,l]}$. Under the smooth-part map, all these different exploded manifolds are sent to  $\{\zeta_{1}\zeta_{2}=0\}$. The discarded  information is parametrized by $\mathbb C^{*}\e{(0,\infty)}$, and is correctly thought of as gluing information.

\subsection{The explosion functor}
Note that $\zeta_{1}\zeta_{2}:\mathbb C^{2}\longrightarrow\mathbb C$ is a normal-crossing degeneration. The explosion functor applied to  a (proper) normal-crossing degeneration  gives a smooth family of exploded manifolds. For a complex manifold $M$ with normal-crossing divisors, $\expl M$ is a holomorphic exploded manifold with smooth part  $M$. The explosion functor replaces coordinate charts that are open subsets $U\subset (\mathbb C^{n},\{\zeta_{1}\dotsb\zeta_{n}=0\})$ with the corresponding open subsets of $\et {n}{[0,\infty)^{n}}$  with smooth part $U\subset \mathbb C^{n}=\totl{\et {}{[0,\infty)^{n}}}$, so we replace coordinates $\zeta_{i}$ (whose vanishing locus give the divisor) with exploded coordinates $z_{i}$ with smooth part $\totl {z_{i}}=\zeta_{i}$. Given any holomorphic map $h:M\longrightarrow N$ sending the interior of $M$ to the interior of $N$ and  each stratum\footnote{The strata of $(\mathbb C^{n},\{\zeta_{1}\dotsc,\zeta_{n}=0\})$ are the sets where $\zeta_{i}=0$ for $i\in I$ and $\zeta_{j}\neq 0$ for $j\notin I$. The strata of $M$ are connected and locally isomorphic to these strata.} of  $M$ into some stratum of $N$, there is a unique explosion of $h$ so that the following diagram commutes.
\[\begin{tikzcd}\expl M\dar\rar{\expl h}&\expl N\dar
\\ \totl{\expl M}=M\rar{h}&\totl{\expl N}=N\end{tikzcd}\]
Deligne-Mumford space has a natural structure of a complex orbifold with normal-crossing divisors the boundary divisors. It turns out that $\expl \bar M_{g,n}$ represents the moduli stack of stable exploded curves with genus $g$ and  $n$ punctures, and the explosion of the forgetful map $\bar M_{g,n+1}\longrightarrow \bar M_{g,n}$ is the universal curve over this moduli stack. As discussed in Remark \ref{log remark}, the explosion functor also applies to nice log schemes (including log smooth log schemes). If $X$ is such a log scheme, we can also apply the explosion functor to the moduli stack of curves in $X$. As discussed in \cite{elc}, the result of this is the moduli stack of holomorphic curves in $\expl X$.

To define Gromov--Witten invariants of $M$ relative  its normal-crossing divisor, we can  use holomorphic curves in $\expl M$. We shall see that the moduli stack of such curves is naturally compact, and (modulo the usual transversality issues) is an orbifold in the category of exploded manifolds. A virtual fundamental class and Gromov--Witten invariants are naturally defined, and are invariant in families of exploded manifolds. In particular, by using the explosion of a normal-crossing  or log smooth degeneration, we get a degeneration formula for Gromov--Witten invariants by calculating Gromov--Witten invariants of an exploded manifold with smooth part the singular fiber.  

The explosion functor works for complex manifolds with normal-crossing divisors, or (nice) log schemes. There is no completely functorial analogue of the explosion functor for symplectic manifolds with normal-crossing divisors, however, by making a contractible choice of extra structure (such as that discussed in \cite{zingerNC}) we can construct an explosion of such a symplectic manifold with normal-crossing divisors. Gromov--Witten invariants of this exploded manifold are invariants of the original symplectic manifold with normal-crossing divisors, so are suitable for defining Gromov Witten invariants relative normal-crossing divisors in the symplectic setting. 
The degeneration formula also works for the symplectic analogue of normal-crossing or log smooth degenerations.

\subsection{Tangent space}

Our defining sheaf of exploded functions $\mathcal E^{\times}(\ex B)$ only has the operation of multiplication, so is not suitable for defining tangent vectors as derivations. Let $\mathcal E(\ex B)$ be the sheaf of $\mathbb C\e{\mathbb R}$--valued functions generated from $\mathcal E^{\times}(\ex B)$ by allowing addition. We can then define vectorfields $v$ as derivations on $\mathcal E(\ex B)$ satisfying the following conditions
\[v(cx+dy)=cvx+dvy\text{ for }x,y\in \mathcal E(\ex B)\text{ and }c,d\in\mathbb C\e{\mathbb R}\ ,\]
\[v(xy)=(vx)y+xvy\ ,\]
\[\totb{v(1\e 0)}=\totb{0\e 0}\ ,\]
and if $x$ is $\mathbb R\e{\mathbb R}$ valued, then so is $vx$.  These conditions ensure that $\totb{vx}=\totb x$ and that if $x\in\mathcal E^{\times}$, then $vx$ is a smooth $\mathbb C$--valued function times $x$.

With these definitions, in standard coordinates $z_{1},\dotsc,z_{m}$ on $\et mP$, the tangent vectorfields are real smooth functions times the real and imaginary parts of $z_{i}\frac{\partial}{\partial z_{i}}$ --- sections of a vectorbundle $T\et mP$ isomorphic to $\mathbb R^{2m}\times\et mP$. Vectorfields and the tangent functor behave largely as they do for smooth manifolds:
\begin{itemize}
\item The relationship between $1$--parameter groups of isomorphisms and flows of vectorfields is as usual.
\item  The cotangent bundle and tensors are defined as usual. Construction of tensors using partitions of unity also work as usual (assuming $\totl {\ex B}$ is second countable).
\item Transversality, intersection, and fiber products of transverse maps work as usual.\footnote{ Transversality and intersection products are determined using $T\ex B$, and not  $\totl{\ex B}$ --- using $\totl{\ex B}$ will not give a correct gluing formula. As a  consequence,   we need  a new cohomology theory, called refined cohomology.}
\end{itemize}

The first surprise about $T\ex B$ is the existence of integral vectors ${}^{\mathbb Z}T_{p}\ex B\subset T_{p}\ex B$. A vector $v$ is integral if $vx/x\in \mathbb Z$ for all exploded functions $x\in\mathcal E^{\times}(\ex B)$. An example of a nonzero integral vector is the real part of $z\frac \partial{\partial z} $ at any point in $\et {}{[0,\infty)}$ with tropical part in $(0,\infty)\subset[0,\infty)$. An arbitrary exploded function, $x=h(\totl{z})\e az^{n}$, has  $vx/x=n$. Over the interior of strata with $m$--dimensional tropical part, ${}^{\mathbb Z}T_{p}\ex B$ is a $\mathbb Z^{m}$--lattice inside $T_{p}\ex B$. We can identify this lattice with the lattice of integral vectors at the image of $p$ in $\totb{\ex B}$, and the obvious identification is functorial.

\subsection{Almost-complex structures and holomorphic curves.}

\begin{defn} An almost-complex structure $J$ on $\ex B$ is an endomorphism of $T\ex B$ (represented as a section of $T^{*}\ex B\otimes T\ex B$) so that $J^{2}=-\id$, and so that, for any integral vector $v$ and exploded function $x$,
\[(Jv)x=i(vx)\]
\end{defn}

For example,  on $\et {}{[0,\infty)}$, $J$ of the real part of $z\frac \partial{\partial z}$ is the imaginary part of $z\frac\partial{\partial z}$ plus a vectorfield  vanishing where $\totl z=0$.

Before we define what a holomorphic curve is, we need the appropriate analogue of `compact' for exploded manifolds:
\begin{defn}$\ex B$ is complete if $\totl{\ex B}$ is compact, and every polytope in $\totb{\ex B}$ is complete (i.e. $\ex B$ is locally isomorphic to $\mathbb R^{n}\times \et {}P$ where $P\subset \mathbb R^{m}$ is closed). 
\end{defn}

For example, $\et {}{(0,1)}$ is not complete, but $\ex T$ is. The explosion of any compact manifold with normal-crossing divisors is also complete. There is a similar notion of a complete map, which is the analogue of a proper map. See Definition 3.5 of \cite{iec}.
\begin{defn}\label{curve} A holomorphic curve is a complete, $2$--dimensional exploded manifold $\ex C$ with an almost-complex structure $j$.

A  curve in $(\ex B,J)$ is a holomorphic curve $(\ex C,j)$ with a map
\[f:\ex C\longrightarrow \ex B\ .\]
This map is a holomorphic curve if $df\circ j=J\circ df$.
\end{defn}

Apart from the exceptional example of $\ex C=\ex T$, each holomorphic curve $(\ex C,j)$ is locally isomorphic to one of three models: 
\begin{itemize}
\item an open subset of $\mathbb C$ --- here our curve behaves like  smooth holomorphic curve; the tropical part of such a connected open subset is just a point;
\item an open subset of $\et {}{[0,\infty)}$ --- the smooth part of our curve here is an open subset of $\mathbb C$ with the special point $0$, the tropical part is $[0,\infty)$, and  we call the strata of $\ex C$ over $(0,\infty)$ an external edge or end of $\ex C$;
\item an open subset of $\et {}{[0,l]}$ --- the smooth part of our curve here is isomorphic to an open subset of $\{\zeta_{1}\zeta_{2}=0\}\subset\mathbb C^{2}$, the model for a node. The strata over this node, $\et {}{(0,l)}$ is called an internal edge of $\ex C$.
\end{itemize}
The smooth part $\totl{\ex C}$ of our curve is some compact nodal curve with a finite collection of special points. 
The tropical part $\totb{\ex C}$ is a complete integral-affine graph with a vertex for every component of $\totl{\ex C}$, an internal edge isomorphic to $[0,l]$ for each node, and a external edge isomorphic to $[0,\infty)$ for each special point.  The extra information in $\ex C$, not seen in the nodal curve $\totl{\ex C}$, is a $\mathbb C^{*}\e{(0,\infty)}$--worth of gluing information for each node, discussed in Example \ref{et[0,l]}.

Because $f:\ex C\longrightarrow \ex B$ has more information than $\totl f:\totl{\ex C}\longrightarrow \totl{\ex B}$, sometimes $\totl f$ has more automorphisms than $f$.\footnote{The cases in which $\totl f$ may have more automorphisms are those in which gluing information at nodes is not preserved. The only case in which $f$ has more automorphisms than $\totl{f}$ is when $\ex C=\ex T$ and $f$ is not injective.} Call $f$ \emph{stable} if both $f$ and $\totl f$ only have a finite number of automorphisms. 

\subsection{Families of exploded manifolds and families of curves}

\begin{defn}A family of exploded manifolds is a complete map $\pi:\hat{\ex B}\longrightarrow \ex B_{0}$ so that for all $p\in \hat{\ex B}$,
\[T_{p}\pi:T_{p}\hat{\ex B}\longrightarrow T_{\pi(p)}\ex B_{0}\text{ is surjective,}\]
\[\text{ and }T_{p}\pi\lrb{{}^{\mathbb Z}T_{p}\hat{\ex B}}={}^{\mathbb Z}T_{\pi(p)}\ex B_{0}\ .\]
\end{defn}
For example, exploding any proper normal-crossing degeneration gives a family of exploded manifolds. 

As usual, there is a notion of vertical (co)tangent bundle on $\hat{\ex B}\longrightarrow \ex B_{0}$, and  a family of tensors  is a section of the appropriate vertical tensor bundle.  Using fiber products, we can perform base changes of families as usual, and  as usual, any family of tensors pulls back under such base changes. As well as holomorphic curves in a single exploded manifold $\ex B$, we are interested in holomorphic curves in a family of exploded manifolds $\hat {\ex B}\longrightarrow \ex B_{0}$ with a family of almost-complex structures $J$.

\begin{defn}\label{curve family} A family of (stable, holomorphic) curves $\hat f$ in a family of targets $(\hat B,J)$ is a family $\ex C(\hat f)\longrightarrow \ex F(\hat f)$ with a family of almost-complex structures $j$, and a map $\hat f$ so that the following diagram commutes
\[\begin{tikzcd}(\ex C(\hat f),j)\rar{\hat f}\dar&(\hat{\ex B},J)\dar
\\ \ex F(\hat f)\rar&\ex B_{0} \end{tikzcd}\] 
and so that $\hat f$ restricted to any fiber is a (stable, holomorphic) curve in the corresponding fiber of $\hat{\ex B}\longrightarrow \ex B_{0}$. 
\end{defn}

\begin{remark}\label{C}The natural regularity for families of curves, and indeed individual curves\footnote{We can achieve more regularity for families of curves if $J$ is integrable, and it may be true that `smooth' is the correct regularity for families of curves when $J$ satisfies an extra condition called being $\dbar$--log compatible. I am unaware of a proof or a counter example.} is a tiny bit weaker than smooth. There is a notion of $\C\infty1$ regularity, defined in section 7 of \cite{iec}, regularity which, for all practical purposes, is as good as smooth. We can still take derivatives to all orders with $\C\infty1$ maps, and $\C\infty1$ exploded manifolds form a category that works  as well as smooth exploded manifolds --- except it takes longer to define. The almost-complex structure $j$ and all maps in Definitions \ref{curve} and \ref{curve family} should be regarded as $\C\infty1$, although for proving compactness of the moduli stack of $J$--holomorphic curves in \cite{cem}, we use stronger regularity on $(\hat {\ex B},J)\longrightarrow \ex B_{0}$.

 On $\et{}{[0,\infty)}$, a $\C\infty 1$ function is a continuous function, smooth where $\totl z\neq 0$, and with  derivatives converging as $\totl z\rightarrow 0$ faster than $\abs{\totl z}^{\delta}$ for every $\delta<1$. Remembering that our vectorfields on $\et {}{[0,\infty)}$ are generated by the real and imaginary parts of $z\frac\partial{\partial z}$, we can consider $\mathbb C\setminus 0$ as having a cylindrical end at $0$, then,  a $\C\infty1$ function is one with all derivatives converging exponentially with all weights $<1$. This $\C\infty1$ is the same as smooth on $\mathbb R$ and $\ex T$, and all other coordinate charts can be obtained as fiber products from $\mathbb R$, $\ex T$, and $\et {}{[0,\infty)}$, and these fiber products are also valid in the category of $\C\infty1$ exploded manifolds.
\end{remark}

\begin{defn} A morphism of $\C\infty1$ families of curves $\hat f\longrightarrow \hat g$  is a commutative diagram of $\C\infty1$ maps
\[\begin{tikzcd}\ex F(\hat f)\dar&\lar (\ex C(\hat f),j)\dar\rar{\hat f}&\hat {\ex B}
\\ \ex F(\hat f)&\lar (\ex C(\hat f),j)\ar{ur}[swap]{\hat g}\end{tikzcd}\]
so that the lefthand square is a pullback diagram.
\end{defn}

\subsection{The moduli stack of curves}

Use the notation $\Msw(\hat {\ex B})$ for the moduli stack of not-necessarily-holomorphic stable curves in $\hat{\ex B}$. This is the category of $\C\infty1$ families of stable curves in $\hat{\ex B}$ along with the functor $\ex F$ assigning the exploded manifold $\ex F(\hat f)$ to the family of curves $\hat f$. Use  $\Mod(\hat {\ex B})\subset \Msw(\hat {\ex B})$ for the substack of holomorphic curves. 

Think of a family of curves parametrized by $\ex F$ as  a map $\ex F\longrightarrow \Msw$. In the case of curves mapping to a point, $\Msw(\cdot)=\Mod(\cdot)$, and both are represented by (the stack of maps into) $\coprod_{g,n}\expl \bar M_{g,n}$.

$\Msw$ comes with a natural topology. Say that a subcategory $\mathcal U$ of $\Msw$ is a substack if $\hat f$ is in $\mathcal U$ if and only if all the individual curves in $\hat f$ are isomorphic to a curve in $\mathcal U$. Define $\mathcal U$ to be open if, for all $\hat f$ in $\Msw$, the subset of $\ex F(\hat f)$ parametrizing curves in $\mathcal U$ is open. For example $\mathcal M\subset \Msw$ is closed, and examples of open substacks are the moduli stack of curves contained in an open subset of $\ex B$, or the stack of curves $f$ so that $\totb f$ has no internal edges.  One interesting example of a closed substack is the substack of $\Msw(\ex B)$ consisting of curves $f$ with at least one component of $\totl{\ex C(f)}$ unstable  and representing a trivial homology class in $\totl{\ex B}$ either by itself or with other components. It's good to know that there is an open neighborhood of $\mathcal M$ excluding those annoying little buggers.

 Interesting phenomena occur in the natural topology on the moduli stack of (not necessarily stable) $\C\infty1$  curves in $\ex B$. The unstable curves are a closed substack, but every neighborhood of an unstable curve includes any stabilization.

\subsection{Regularity of the moduli stack of holomorphic curves}

\

There is a reasonably nice (but infinite-dimensional) bundle $\Y$ over $\Msw(\hat {\ex B})$ discussed in Section 2.2 of \cite{evc}. This bundle comes with a section $\dbar:\Msw\longrightarrow \Y$ defining holomorphic curves:  $\Mod\subset \Msw$ is the intersection of $\dbar$ with $0$. If $\dbar$ is transverse\footnote{See \cite{evc}, sections 2.7 and 2.8 for a discussion of the tangent space of $\Msw$ and the linearization of $\dbar$.} to $0$ at a holomorphic curve $f$, then there exists an open neighborhood $\mathcal U\subset \Msw$ of $f$, and  a $\C\infty1$ family of curves, $\hat f$, with group of automorphisms, $G$, so that $\Mod\cap \mathcal U$ is represented by $\hat f/G$. So, the moduli stack of curves is an orbifold close to where $\dbar$ is transverse to $0$.

In general, $\dbar$ is not transverse to $0$. However,  on a small-enough neighborhood  $\mathcal U$ of $f$, there always exists a suitably nice\footnote{Definition 2.24 of \cite{evc}} finite-dimensional vector sub-bundle $V$ of $\Y$ so that $\dbar$ is transverse to $V$. Theorem 6.6 of \cite{evc} then states that $\dbar^{-1} V \subset \mathcal U$ is an orbifold locally represented by $\hat f/G$ for some family of curves $\hat f$. Then $V(\hat f)$ is a finite-dimensional, $G$--equivariant vectorbundle over $\ex F(\hat f)$ with a natural section 
\[\dbar\hat f:\ex F(\hat f)\longrightarrow V(\hat f)\]
and the moduli stack of holomorphic curves is locally represented by the quotient by $G$ of the intersection of $\dbar\hat f$ with $0$.
The data $(\mathcal U,V,\hat f/G)$ with the section $\dbar$ defines a  Kuranishi chart for $\Mod\subset\Msw$, a Kuranishi chart embedded in $\Msw$. We can also choose $V$ to be a complex vectorbundle, and $\dbar$ to be transverse to $0$ in a stronger sense\footnote{Definition 2.26 of \cite{evc}} so that the homotopy of $D\dbar$ to a $\mathbb C$--linear operator defines a canonical orientation of $\ex F(\hat f)$ (relative to $\ex B_{0}$ in the case of curves in a family of targets $\hat{\ex B}\longrightarrow \ex B_{0}$).
\begin{remark}The statement of regularity above includes all necessary gluing statements, as a neighborhood of a nodal curve includes all relevant gluings.  It also includes the gluing statements required for relating the Gromov--Witten invariants of a fiber  of $\hat{\ex B}\longrightarrow \ex B_{0}$  to the invariants of nearby fibers. 
\end{remark}

Under the assumption that $\Mod\subset\Msw$ obeys suitable compactness conditions, Theorem 7.3 of \cite{evc} constructs a compatible collection\footnote{We are omitting many details. Roughly speaking, we construct a collection of embedded Kuranishi charts, compatible in the sense that on their common domain of definition, one includes into the other. One property of our Kuranishi charts we have not mentioned is that they can be shrunk appropriately. This property is essential for making global constructions of sections of sheaves over Kuranishi charts, and is related to the Hausdorff condition one could desire for making such global constructions. } of such Kuranishi charts covering $\Mod$. In \cite{vfc}, a virtual fundamental class is constructed from such a Kuranishi structure.

\subsection{Taming forms and compactness of the moduli stack of holomorphic curves}

\

The smooth part $\totl{\ex B}$ of $\ex B$ also has a (co)tangent space, with a canonical complex structure normal to strata.\footnote{See section 2 of \cite{cem}. Unlike $T\ex B$, the tangent space of $\totl{\ex B}$ is not always a vectorbundle, and is not particularly nice.} 

A taming form $\omega$ on $\ex B$ is a closed $2$--form on $\totl{\ex B}$, symplectic on strata, and positive on holomorphic planes normal to strata.\footnote{For a more precise statement, see section 2 of \cite{cem}.}

An almost-complex structure on $\ex B$ does not necessarily induce an almost-complex structure on $\totl{\ex B}$. A nice kind of almost-complex structure that does induce an almost-complex structure on $\totl{\ex B}$ is called a $\dbar$--log compatible complex structure. We can speak of $\omega$ taming such a $\dbar$-compatible $J$ if $\omega$ is positive on  holomorphic planes within $T\totl{\ex B}$. Roughly speaking, $J$ is $\dbar$--log compatible if for all exploded functions $x\in\mathcal E^{\times}\ex B$, 
\[f^{-1}(df-idf\circ J)\]
is a $\mathbb C$--valued $1$--form on $\ex B$ pulled back from $\totl{\ex B}$.\footnote{ See section 3 of \cite{cem}. Note that any smooth almost-complex structure on a smooth manifold is also   $\dbar$--log compatible. } The standard complex structure on $\et{}P$ is $\dbar$--log compatible, and the space of $\dbar$--log compatible almost-complex structures on $\ex B$ tamed by a given taming form is nonempty and connected. 

The following theorem follows from Theorem 6.1 and Lemma 4.2 of \cite{cem} along with Lemma 2.8 of \cite{evc}. 

\begin{thm}Suppose that $J$ is a $\dbar$--log compatible almost-complex structure tamed by $\omega$ on a complete exploded manifold $\ex B$, and suppose further that there exists an affine immersion $\totb{\ex B}\longrightarrow [0,\infty)^{N}$. Then the moduli stack of holomorphic curves in $\Msw(\ex B)$ with genus $g$, $n$ ends, and $\omega$--energy at most $E$ is compact.\footnote{The reader may wonder why we use `compact' here instead `complete'. The virtual moduli space is complete, but $\Mod$ may not be complete when $\dbar$ is not transverse to $0$. When we take the non-transverse intersection of the zero-section with  a section of a vectorbundle over a complete exploded manifold, it may not be complete, but it is always compact.}

If there only exists an immersion $\totb{\ex B}\longrightarrow \mathbb R$, then $\mathcal M$ is compact restricted to every connected component of $\Msw(\ex B)$.
\end{thm}
The analogous result holds for holomorphic curves in  a family of exploded manifolds $\hat{\ex B}\longrightarrow \ex B_{0}$, but now $\mathcal M\longrightarrow\ex B_{0}$ is proper when restricted to curves with bounded energy, number of ends and genus or when restricted to connected components of $\Msw(\hat {\ex B})$ respectively. 

\subsection{DeRham cohomology theories}

\

There are three different DeRham cohomology theories\footnote{$H^{*}$ and $\rh^{*}$ are defined in \cite{dre}. $\rhf^{*}$ is defined in \cite{vfc}.} for exploded manifolds: ordinary cohomology $H^{*}$, refined cohomology $\rh^{*}$ and refined cohomology generated by functions, $\rhf ^{*}$. All  are identical to usual DeRham cohomology when applied to smooth manifolds, and have integration and Stokes' theorem working as usual. These $3$ cohomology theories are related as follows:
\[H^{*}(\ex B)\hookrightarrow \rh^{*}(\ex B)\longleftarrow \rhf^{*}(\ex B)\]
Ordinary cohomology $H^{*}(\ex B)$ is finite-dimensional if $\ex B$ is compact, invariant\footnote{More accurately, $H^{*}$ of the fibers of a family of exploded manifolds forms a kind of flat bundle over the family with interesting monodromy. See section 11 of \cite{dre}.} in connected families of exploded manifolds, and for $M$ a complex manifold with normal-crossing divisors, $H^{*}(\expl M)=H^{*}(M,\mathbb R)$. The problem with ordinary cohomology is that it does not have pushforwards compatible  with  fiber products of exploded manifolds, so it is inappropriate for our gluing formula. 
  
 Refined cohomology,  $\rh^{*}$ is the minimal extension of $H^{*}$ with fiber-product-compatible pushforwards.  It is also invariant in connected families of exploded manifolds and is usual cohomology on smooth manifolds, but when $M$ is a complex manifold with normal-crossing divisors, $\rh^{*}(\expl M)$ is usually infinite dimensional. 
 
 Refined cohomology generated by functions $\rhf^{*}$ also possesses fiber-product-compatible pushforwards, but it is only invariant in families of exploded manifolds parametrized by $\mathbb R$, not general families. This non-invariance is actually a convenient feature, because it allows multiple different gluing formulas for the same Gromov--Witten invariant. The main advantage of $\rhf^{*}$ is its compatablity with tropical completion, discussed in section 7 of \cite{vfc}. For a given point $p\in\totb{\ex B}$, the set of points in $\ex B$ with tropical part $p$ is a  manifold, $\ex B\rvert_{p}$, with a canonical completion, $\ex B\tc p$, to an exploded manifold that is complete when $\ex B$ is compact.  Any differential form, $\theta$ on $\ex B$, restricts to a differential form on $\ex B\rvert _{p}$ and then extends canonically to a differential form $\theta\tc p$ on $\ex B\tc p$. The map $\theta\mapsto \theta\tc p$ determines a well-defined map
  \[\rhf^{*}(\ex B)\longrightarrow \rh^{*}(\ex B\tc p)\ .\]
 
 \subsection{Gromov--Witten invariants}
 
 \
 
 A virtual fundamental class $[\Mod]$ for the moduli stack of holomorphic curves is constructed in \cite{evc}. We can integrate forms from $H^{*}$, $\rh^{*}$ or $\rhf^{*}$ over $[\Mod]$ to define numerical Gromov--Witten invariants. Given any evaluation map
 \[ev:\Msw(\ex B)\longrightarrow \ex X\]
 where $\ex X$ is an oriented exploded orbifold or manifold and $ev$ is proper restricted to $\Mod_{g,n,E}$ --- the moduli stack of holomorphic curves with genus $g$, $n$ ends, and $\omega$--energy $E$ --- we can pull back and integrate $\theta\in \rh^{*}_{c}(\ex X)$.\footnote{We can also use differential forms representing classes in $H^{*}$ or $\rhf^{*}$, as these differential forms automatically represent classes in $\rh^{*}$.} The integral 
 \[\int_{[\Mod_{g,n,E}]}ev^{*}\theta\]
 is well defined independent of all choices, and invariant in families $\hat{\ex B}\longrightarrow \ex B_{0}$ (so long as the evaluation map makes sense for these families).  We can also push forward the virtual fundamental class to a class $\eta_{g,n,E}\in \rhf^{*}(\ex X)$ so that
 \[\int_{[\Mod_{g,n,E}]}ev^{*}\theta=\int_{\ex X}\eta_{g,n,E}\wedge \theta\ .\]
 This class $\eta_{g,n,E}\in \rhf^{*}(\ex X)$ is independent of all choices and invariant in families parametrized by $\mathbb R$. Its image in $\rh^{*}(\ex X)$ is invariant in connected families in the following sense. Given a family of evaluation maps to a family of targets, relatively oriented over $\ex B_{0}$,
 \[\begin{tikzcd}\Msw(\hat{\ex B})\dar\rar{ev}&\hat{\ex X}\dar
 \\ \ex B_{0}\rar{\id}&\ex B_{0}
 \end{tikzcd}\]
 so long as $ev$ is proper restricted to $\Msw_{g,n,E}(\hat{\ex B})$, we can again push forward the virtual fundamental class to define
  $\eta_{g,n,E}\in \rhf^{*}(\hat{\ex X})$.\footnote{For tagging the moduli space with $\omega$--energy $E$, we are assuming the the cohomology class represented by the taming form $\omega$ does not change in our family $\hat{\ex B}\longrightarrow \ex B_{0}$. An analogous statement holds if this is not the case. } The precise statement of invariance is that this class, $\eta_{g,n,E}$, is compatible with base changes. This base-change compatibility follows from the observation that embedded Kuranishi structures pull back under base changes, and Theorem 5.22 of \cite{evc}. 
 
 For stating our gluing theorem, we package the different $\eta_{g,n,E}$ into a generating function
 \[\eta:=\sum_{g,n,E}\hbar^{2g+n-2}q^{E}\eta_{g,n,E}\]
 \begin{remark}It is also possible to encapsulate the pushforward of tautological classes in $\eta$ so that it also records descendant Gromov--Witten invariants. The same invariance properties still hold, and the gluing formula still works. See \cite{tpgf} for an example of this. \end{remark}
 \subsection{Tropical gluing formula}
 
 \
 
 Let $\gamma$ be a tropical curve in $\totb{\ex B}$, and let $\Msw_{\gamma}(\ex B)$ denote the moduli stack of stable curves in $\ex B$ with a chosen isomorphism of their tropical part with $\gamma$. This stack is an $\abs{\Aut\gamma}$--fold cover of the stack of curves in $\Msw(\ex B)$  with tropical part isomorphic to $\gamma$. The contribution of $\Msw_{\gamma}$ to $\eta$ is well defined:
 \[\text{ for }p\in\totb{\ex X},\ \ \  \eta\tc p=\sum_{\gamma\in\totb{ev}^{-1}(p)}\eta_{\gamma}/\abs{\Aut\gamma}\]
 And, $\Msw_{\gamma}$ is a fiber product of moduli stacks explained below.
\begin{equation}\label{Bcd}\begin{tikzcd}[column sep=large]\Msw_{\gamma}(\ex B)\rar \dar{\text{cut}}&\prod_{e\in\ie(\gamma)}\ex B\rvert_{e}/\mathbb C^{*}\dar{\Delta}
\\ \prod_{v}\Msw_{\gamma_{v}}(\ex B\tc v)\rar&\prod_{e\in\ie(\gamma)}(\ex B\rvert_{e}/\mathbb C^{*})^{2}\end{tikzcd}\end{equation}
For any vertex $v$ of $\gamma$, we can obtain a curve in $\ex B\tc v$ by applying tropical completion at $v$ to a curve in $\Msw_{\gamma}(\ex B)$. The subset of the domain $\ex C$ over $v$ is a punctured Riemann surface $\ex C\rvert_{v}$, which maps to $\ex B\rvert_{v}$. This map has a canonical completion to a curve  $\ex C\tc v\longrightarrow \ex B\tc v$, with tropical part parametrized by $\gamma_{v}$ --- the tropical curve obtained from $\gamma$ by cutting at all internal edges, then semi-infinitely extending the cut edges of the component containing $v$. Applied to a family of curves in $\Msw_{\gamma}(\ex B)$, this process creates a family of curves in $\Msw_{\gamma_{v}}(\ex B\tc v)$ for all $v$, and determines a map of stacks $\text{cut}:\Msw_{\gamma}(\ex B)\longrightarrow \prod \Msw_{\gamma_{v}}(\ex B\tc v)$. 

 In the above diagram, $\ie(\gamma)$ indicates the set of internal edges of $\gamma$. Over each point in an edge of $\gamma$, the domain of a curve $f$ has  a copy of $\mathbb C^{*}$.  Similarly, over every point $p$ in a $n$-dimensional stratum of $\totb{\ex B}$, the manifold $\ex B\rvert_{p}$,  is a $(\mathbb C^{*})^{n}$--bundle. Pick a point in the tropical part of each edge $e$ of $\gamma$, and let $\ex B\rvert_{e}\subset \ex B$ be the subset of $\ex B$ with tropical part equal to this point.  The restriction of $f$ determines a map $\mathbb C^{*}\longrightarrow \ex B\rvert_{e}$  equivariant with respect to a $\mathbb C^{*}$--action on $\ex B\rvert_{e}$ --- the $\mathbb C^{*}$--action  with weight $(a_{1},\dotsc,a_{n})$ when the derivative of $\gamma$ on $e$ is $(a_{1},\dotsc,a_{n})$. The top line of Diagram \ref{Bcd} is the associated  evaluation map. 
\begin{equation}\label{7}\Msw_{\gamma}(\ex B)\longrightarrow \prod_{e\in \ie (\gamma)}\ex B\rvert_{e}/\mathbb C^{*}\end{equation}
There is a similar evaluation map for each edge of $\gamma_{v}$. For each internal edge of $\gamma$,  there are two corresponding edges  in $\coprod_{v}\gamma_{v}$, so  the bottom line of Diagram \ref{Bcd} has target $\prod (\ex B\rvert_{e}/\mathbb C^{*})^{2}$.

In each case $\ex B\rvert_{e}$ is the quotient of a manifold $X_{e}$ by a trivial group action.
Use $m_{e}$ to denote the number so that derivative of $\gamma$ at the edge $e$ is $m_{e}$ times a primitive integral vector. If $m_{e}$ is zero, the $\mathbb C^{*}$--action on $\ex B\rvert_{e}$ is trivial, and  $X_{e}:=\ex B\rvert_{e}$. Otherwise, $\ex B\rvert_{e}/\mathbb C^{*}$ is the quotient of a manifold $X_{e}$ by the trivial $\mathbb Z_{m_{e}}$ action.  

\begin{remark}Each of the stacks in Diagram \ref{Bcd} is actually a stack over the category of manifolds. The associated moduli stacks of holomorphic curves are not compact. We really need to work with appropriate completions of these stacks in order to prove our gluing theorem; see \cite{gfgw}.   \end{remark}

$\Msw_{\gamma_{v}}(\ex B\tc v)$ has a natural closure $\Msw_{[\gamma_{v}]}(\ex B\tc v)$ in the moduli stack of curves with labeled ends, and $X_{e}$ sits inside an exploded manifold $\rend (\ex B\tc v)$ with a similar evaluation map.
\[\begin{tikzcd}\Msw_{\gamma_{v}}(\ex B\tc v)\dar[hook]\rar& X_{e}\dar[hook]
\\ \Msw_{[\gamma_{v}]}(\ex B)\rar& \rend(\ex B\tc v)\end{tikzcd}\]
The corresponding evaluation map at all $k$ edges of $\gamma_{v}$ gives a map to a component, $\rend_{\gamma_{v}}(\ex B\tc v)$, of $(\rend(\ex B\tc v))^{k}$. There is a similar evaluation map from $\Msw (\ex B)$, recording the location of all external edges. As we don't have a labeling of external edges, this evaluation map takes values in $\coprod_{n}(\rend\ex B)^{n}/S_{n}$. We now have the elements of the key diagram for our gluing formula.
\[\begin{tikzcd}\Msw(\ex B)\rar{ev}&\coprod_{n}(\rend \ex B)^{n}/S_{n}
\\ \Msw_{\gamma}\uar\rar\dar{\text{cut}}&\prod_{\ed(\gamma)} X_{e}\dar{\Delta}\uar{i_{\gamma}}
\\ \prod_{v}\Msw_{[\gamma_{v}]}(\ex B\tc v)\rar{\prod_{v}ev_{\gamma_{v}}}&\prod_{v}\rend_{\gamma_{v}}(\ex B\tc v)\end{tikzcd}\]
In the above,  $\ed(\gamma)$ indicates all edges of $\gamma$, and $i_{\gamma}$ forgets $X_{e}$ for internal edges and otherwise uses the inclusion of $X_{e}\subset \rend\ex B$ for external edges.
 
 The gluing formula reads\footnote{We have cheated a little bit. To know that the pushforward $(i_{\gamma})_{!}$ gives a well-defined cohomology class, we need an extension to a map from a compact exploded manifold containing $X_{e}$. There are some complications involving tropical completions. Further details are given in \cite{gfgw}.} 
 \[\eta \tc p=\sum_{\gamma\in \totb{ev}^{-1}(p)}\frac {k_{\gamma}}{\abs{\Aut\gamma}}(i_{\gamma})_{!}\Delta^{*}\bigwedge \eta_{\gamma_{v}}\]
 where $\eta_{\gamma_{v}}\in\rhf^{*}(\rend_{\gamma_{v}}(\ex B))$ is the Gromov--Witten invariant associated with $ev_{\gamma_{v}}$, and 
 \[k_{\gamma}=\prod_{\ie(\gamma)}m_{e}\]
 is an extra combinatorial factor introduced  after we forgot the trivial group action on $X_{e}$ present in Diagram \ref{Bcd}.
 
 \subsection{Refinements}
 
 \
 
Refinements, and refined cohomology,  are important for understanding the tropical gluing formula in terms of manifolds with normal-crossing divisors, or log schemes instead of exploded manifolds.
 
 \begin{defn} A refinement $\ex B'\longrightarrow \ex B$ is a complete, bijective submersion.
 \end{defn}
 As is explained in section 10 of \cite{iec}, refinements of $\ex B$ are equivalent to subdivisions of the tropical part $\totb{\ex B}$ of $\ex B$. An instructive example is a refinement of $\ex T^{n}$ determined by subdividing $\mathbb R^{n}=\totb{\ex T^{n}}$ into the toric fan of some toric manifold $M$. This refinement is $\expl M$, the explosion of $M$ relative to its toric boundary divisors. A further subdivision\footnote{A conical subdivision of the toric fan is the explosion of the toric manifold with that subdivided fan. A more general subdivision would result in an exploded manifold $\ex B'$ which is not the explosion of a log scheme or complex manifold with normal-crossing divisors. The smooth part $\totl{\ex B'}$ of $\ex B'$ would have extra components. Special cases of these extra components have been variously called  holomorphic buildings, rubber components and expansions.} of this toric fan would result in a further refinement of $\expl M$, which in turn would be the explosion of some toric blowup of $M$. 
 
Refinements do not affect Gromov Witten invariants, suitably interpreted.  (The virtual fundamental class of) the moduli stack of holomorphic curves in $\ex B'$ is a refinement of (the virtual fundamental class of) the moduli stack of holomorphic curves in $\ex B$.\footnote{As well as my proof of this using an older construction of the virtual fundamental class in \cite{egw}, Abramovich and Wise have proved this in the log setting in \cite{ilgw}.}

Refinements are almost isomorphisms. One manifestation  is that the following  is a fiber-product diagram.
 \[\begin{tikzcd}\ex B'\rar\dar&\ex B'\dar
 \\ \ex B'\rar&\ex B\end{tikzcd}\]
So,  any cohomology theory with pushforwards compatible with fiber product diagrams and containing $H^{*}$ must also contain $H^{*}(\ex B')$ for any refinement $\ex B'$. Refined cohomology $\rh^{*}(\ex B)$ is the minimal such cohomology theory. Refinements of $\expl M$ correspond to boundary blowups  $M'$ of $M$ locally modeled on toric blowups. So, refined cohomology\footnote{defined differently, but equivalently in \cite{dre}} of $\expl M$ is the direct limit of the cohomology of all such $M'$.

We can use refinements to identify the key elements of the tropical gluing formula in more familiar terms. First, let us describe the evaluation map to $\rend \ex B$. Suppose that $\ex B$ is the explosion of a complex manifold $M$ relative to a simple normal-crossing divisor $\bigcup_{i=1}^{k} D_{i}$ with irreducible components. And suppose, for simplicity, that the intersection of any number of these $D_{i}$ is connected. (This can be achieved by suitably blowing up $M$, corresponding to a refinement of $\ex B$.)  The tropical part of $\ex B$ may be identified with the vectors $(v_{1},\dotsc,v_{k})\in[0,\infty)^{k}$ so that if $v_{i}>0$ for some set of indices $i\in I$, then $\bigcap_{i\in I}D_{i}\neq \emptyset$.  The connected components of $\rend\ex B$ are then indexed by nonnegative integral vectors in $\mathbb N^{k}$ satisfying the same condition. We can consider these integral vectors as either specifying the derivative at an end of a tropical curve, or  specifying the degree of  contact with each $D_{i}$. The zero vector corresponds to no contact with any $D_{i}$. The corresponding component of $\rend\ex B$ is $\ex B=\expl M$. The vector $(n,0,\dotsc,0)$ corresponds to order $n$ contact with $D_{1}$, and no contact with any other $D_{i}$. The corresponding component of $\rend\ex B$ is $\expl D_{1}$, where we use the divisor given by intersection with the other $D_{i}$. Similarly, the component of $\rend\ex B$ corresponding to order $n$ contact with $D_{i}$ and no contact with other divisors is $\expl D_{i}$. No other component of $\rend\ex B$ will be the explosion of anything, but a refinement of it will be; simply blow up $M$ (refine $\ex B$ to $\ex B'$) until the specified order of contact is with only one component $D$ of the divisor. Then  $\rend\ex B'$ is a refinement of $\rend{\ex B}$, and the refinement of our component is $\expl D$. As $\Mod(\ex B')$ is also a refinement of $\Mod(\ex B)$, we can use the map $\Mod(\ex B')\longrightarrow \expl D\subset \rend{\ex B'}$ to define GW invariants in place of  $\Mod(\ex B)\longrightarrow \rend\ex B$.
 
 The exploded manifolds $\ex B\tc v$ may also be described with the help of refinements. Suppose that $\ex B$ is a fiber of the explosion of a simple normal-crossing degeneration with singular fiber $\cup_{i}N_{i}$. When $\ex B$ is a fiber with nontrivial tropical part, its smooth part is equal to this singular fiber. The tropical part of $\ex B$ has a vertex $p_{i}$ for each $N_{i}$, and a simplex with corners $\{p_{i}, \  i\in I\}$ for each connected component of $\bigcap_{i\in I} N_{i}$. If $v$ is one of these $p_{i}$, then $\ex B\tc v$ is $\expl N_{i}$. When $v$ is in a higher-dimensional stratum of $\totb{\ex B}$, we can refine $\ex B$ to $\ex B'$, subdividing $\totb{\ex B}$ so that $v$ becomes a zero-dimensional stratum.\footnote{We can understand refinements of $\ex B$ in terms of refinements of $\expl M$, where $M$ is the total space of the normal-crossing degeneration.  $\totb{\expl M}$ is a cone over $\totb{\ex B}$. Our refinement is equivalent to a subdivision of $\totb{\ex B}$; if the cone over this subdivision has rational slopes, it is a subdivision defining a refinement of $\totb{\expl M}$. Beware that refining the total space over a family will generally not produce a family because it will no longer satisfy the requirement on integral vectors. This is dealt with by a base change.} So long as every stratum in $\totb{\ex B'}$ with $v$ as a corner has $v$ as a standard corner,\footnote{A standard corner of an integral-affine polytope is locally isomorphic to an open neighborhood of $0$ in  $[0,\infty)^{n}$.} the stratum of $\totl{\ex B}$ corresponding to  $v$ is  a manifold $N$ with normal-crossing divisors. Then  $ \ex B'\tc v=\expl N$ is a refinement of $\ex B\tc v$, and we can use $\ex B'\tc v$ in place of $\ex B\tc v$ for computing Gromov--Witten invariants. 
 
 With the above translations using refinements or blowups, our tropical gluing formula translates into the language of log schemes or ordinary manifolds. 
\bibliographystyle{plain}
\bibliography{ref}

\begin{thebibliography}{10}

\bibitem{acgw}
Dan Abramovich and Qile Chen.
\newblock Stable logarithmic maps to {D}eligne-{F}altings pairs {II}.
\newblock {\em The Asian Journal of Mathematics}, 18(3):465--488, 2014.

\bibitem{ilgw}
Dan Abramovich and Jonathan Wise.
\newblock Invariance in logarithmic {G}romov-{W}itten theory.
\newblock \href{http://arxiv.org/abs/1306.1222}{arxiv:1306.1222}, 2013.

\bibitem{Chen}
Qile Chen.
\newblock Stable logarithmic maps to {D}eligne-{F}altings pairs {I}.
\newblock {\em Ann. of Math. (2)}, 180(2):455--521, 2014.

\bibitem{GSlogGW}
Mark Gross and Bernd Siebert.
\newblock Logarithmic {G}romov-{W}itten invariants.
\newblock {\em J. Amer. Math. Soc.}, 26(2):451--510, 2013.

\bibitem{IonelGW}
Eleny-Nicoleta Ionel.
\newblock G{W} invariants relative to normal crossing divisors.
\newblock {\em Adv. Math.}, 281:40--141, 2015.

\bibitem{zingerNC}
Mark McLean, Mohammad Tehrani, and Aleksey Zinger.
\newblock Normal crossings divisors and configurations for symplectic topology.
\newblock arXiv:1410.0609.

\bibitem{tpgf}
Brett Parker.
\newblock Gluing formula for {G}romov-{W}itten invariants in a triple product.
\newblock \href{http://arxiv.org/abs/1511.0779}{arXiv:1511.0779}.

\bibitem{dre}
Brett Parker.
\newblock De {R}ham theory of exploded manifolds.
\newblock \href{http://arxiv.org/abs/1003.1977}{arXiv:1003.1977}, 2011.

\bibitem{egw}
Brett Parker.
\newblock Gromov-{W}itten invariants of exploded manifolds.
\newblock \href{http://arxiv.org/abs/1102.0158}{arXiv:1102.0158}, 2011.

\bibitem{iec}
Brett Parker.
\newblock Exploded manifolds.
\newblock {\em Adv. Math.}, 229:3256--3319, 2012.
\newblock \href{http://arxiv.org/abs/0910.4201}{arXiv:0910.4201}.

\bibitem{elc}
Brett Parker.
\newblock Log geometry and exploded manifolds.
\newblock {\em Abh. Math. Sem. Hamburg}, 82:43--81, 2012.
\newblock \href{http://arxiv.org/abs/1108.3713}{arxiv:1108.3713}.

\bibitem{evc}
Brett Parker.
\newblock Holomorphic curves in exploded manifolds: Kuranishi structure.
\newblock \href{http://arxiv.org/abs/1301.4748}{arXiv:1301.4748}, 2013.

\bibitem{tropicalIonel}
Brett Parker.
\newblock On the value of thinking tropically to understand {I}onel's {GW}
  invariants relative normal crossing divisors.
\newblock \href{http://arxiv.org/abs/1407.3020}{arXiv:1407.3020}, 2014.

\bibitem{cem}
Brett Parker.
\newblock Holomorphic curves in exploded manifolds: compactness.
\newblock {\em Adv. Math.}, 283:377--457, 2015.
\newblock \href{http://arxiv.org/abs/0911.2241}{arXiv:0911.2241}.

\bibitem{vfc}
Brett Parker.
\newblock Holomorphic curves in exploded manifolds: virtual fundamental class.
\newblock \href{http://arxiv.org/abs/1512.05823}{arXiv:1512.05823}, 2015.

\bibitem{gfgw}
Brett Parker.
\newblock Tropical gluing formulae for {G}romov-{W}itten invariants.
\newblock \href{http://arxiv.org/abs/1703.05433}{arXiv:1703.05433}, 2017.

\end{thebibliography}
 \end{document}